\documentclass[english,a4paper]{article}
\usepackage{babel}
\usepackage{epsfig,graphicx,times}
\usepackage{amssymb,amsmath}

\def\Ps{\mathcal{P}}

\newcommand{\di}{m}

\def\@begintheorem#1#2{\list{}{\thm@body}%
  \item[]{\bf #1~#2.}\quad\it\ignorespaces}
\def\@opargbegintheorem#1#2#3{\list{}{\thm@body}%
  \item[]{\bf #1~#2~\ifrembrks #3\global\rembrksfalse\else (#3)\fi.}%
  \quad\it\ignorespaces}
\def\@endtheorem{\endlist}
\newtheorem{theorem}{Theorem}

\newtheorem{definition}{Definition}

\newcommand{\eop}{\hfill{$\Box$}}
\newcommand{\I}{\mathcal{I}}

\newcommand{\Z}{\mathbb{Z}_n^\di}

\begin{document}

\title{There are integral heptagons, no three points on a line, no four on a circle}
\author{{\sc Tobias Kreisel}\thanks{tobias.kreisel@uni-bayreuth.de} { and }{\sc Sascha Kurz}\thanks{sascha.kurz@uni-bayreuth.de}\\ 
      Department of Mathematics, University of Bayreuth\\ 
      D-95440 Bayreuth, Germany}

\maketitle

\vspace*{-4mm}
\noindent
{\center\small{Keywords: integral distances, exhaustive search, orderly generation, solution to an Erd\H{o}s problem\hspace*{2cm} MSC: 52C10,52C35,52-04,52A99,51K99\\}}
\noindent
\rule{\textwidth}{0.3 mm}

\begin{abstract}
  \noindent
  We give two configurations of seven points in the plane, no three points in a line, no four points on a circle 
  with pairwise integral distances. This answers a famous question of Paul Erd\H{o}s.
\end{abstract}
\noindent
\rule{\textwidth}{0.3 mm}

\section{Introduction}

A famous open problem of P. Erd\H{o}s asks for seven points in the plane, no three on a line, no four on a circle with pairwise rational or integral distances \cite{1086.52001,UPIN}. For six points parameter solutions for infinite families of such point sets are known, see e.g. \cite{hab_kemnitz}. Since for finite point sets we can multiply the occurring distances with their denominators' smallest common multiple we confine ourselves to considering integral distances only. From the combinatorial point of view the question for the smallest possible diameter $\dot{d}(2,n)$ of $n$ points arises, where the diameter is the largest occurring distance in a point set. So far 
$$
  \left(\dot{d}(2,n)\right)_{n=3,\dots,6}=1,8,73,174
$$
are known \cite{integral_distances_in_point_sets}. By exhaustive search the bound $\dot{d}(2,7)\ge 20000$ could be determined \cite{1088.52011,paper_alfred}. Up to diameter $20000$ there are only few integral point sets consisting of $6$ points, no three on a line, no four on a circle with pairwise integral distances, see \cite{hp} for a complete list. Some attempts to show that no integral point set in general position consisting of more than six points can exist are known \cite{pers}, but the suggested proofs turned out to be incorrect. So there was little hope to discover such a point set. But then by a suggestion of S. Dimiev \cite{Dimiev-Setting} we considered integral point sets over $\mathbb{Z}_n^2$ \cite{paper_axel}.
 
\begin{definition}
  Two points $(u_1,\dots,u_\di),(v_1,\dots,v_\di)\in \Z:=(\mathbb{Z}\backslash\mathbb{Z}n)^m$ are at \textbf{integral distance} if 
  there exists a number $d\in\mathbb{Z}_n$ with
  $
    \sum\limits_{i=1}^{\di}(u_i-v_i)^2=d^2
  $.
\end{definition}

\noindent
So, an integral point set in $\mathbb{Z}_n^2$ is defined as a subset of $\mathbb{Z}_n^2$ where all pairs of points are at integral distance. To have an analogue to the ``no three on a line and no four on a circle" restriction we need two further definitions.

\begin{definition}
  \label{definition_collinear}
  A set of $r$ points $(u_i,v_i)\in\mathbb{Z}_n^2$ is collinear if there are $a,b,t_1,t_2,w_i\in \mathbb{Z}_n$ 
  with
  $
    a+w_it_1=u_i\,\,\text{and}\,\, b+w_it_2=v_i
  $.
\end{definition}

\begin{definition}
  Four points $p_i=(x_i,y_i)$ in $\mathbb{Z}_n^2$ are said to be situated on a circle if there exist $a,b \in\mathbb{Z}_n$, 
  $r \in\mathbb{Z}_n\backslash\{\overline{0}\}$ with 
  $
    (x_i-a)^2+(y_i-b)^2=r^2\,\,\forall i
  $.
\end{definition}

\noindent
By $\dot{\I}(n,2)$ we denote the maximum number of points in $\mathbb{Z}_n^2$ with pairwise integral distances where no 
three are collinear and no four points are situated on a circle. By combinatorial search techniques---see \cite{paper_axel} 
for the details---we found two point sets proving $\dot{\I}(50,2)\ge 12$ and $\dot{\I}(61,2)\ge 9$. Surely this does not 
imply the existence of an integral point set over the real plane in general position, i.e. no three points on a line, no four points on a 
circle, however it did give us a fresh impetus to continue our search.

\section{Integral heptagons in general position}

The results for the ``relaxed" problem over $\mathbb{Z}_n^2$ motivated us to maintain our approach of exhaustive generation of all plane integral point sets in general position up to a given diameter by a variant of orderly generation, see \cite{1088.52011,paper_alfred} for details. Also, without changing our approach but simply by harnessing more computational power we were lucky enough to discover the following distance matrix
\begin{equation}
  \label{example_1}
  \left(
  \begin{array}{rrrrrrr}
        0 & 22270 & 22098 & 16637 &  9248 &  8908 &  8636 \\
      22270 &     0 & 21488 & 11397 & 15138 & 20698 & 13746 \\
      22098 & 21488 &     0 & 10795 & 14450 & 13430 & 20066 \\
    16637 & 11397 & 10795 &     0 &  7395 & 11135 & 11049 \\
     9248 & 15138 & 14450 &  7395 &     0 &  5780 &  5916 \\
     8908 & 20698 & 13430 & 11135 &  5780 &     0 & 10744 \\
     8636 & 13746 & 20066 & 11049 &  5916 & 10744 &     0
  \end{array}
  \right)
\end{equation}

\noindent
corresponding to a plane integral point set in general position with diameter $22270$ consisting of seven points. So this 
answers Erd\H{o}s's question positively. Since we applied an exhaustive search we receive:

\begin{theorem}
  $
    \dot{d}(2,7)=22270
  $.
\end{theorem}

\noindent
To avoid duplicated listings of isomorphic point sets we give all point sets in the following canonical form. Consider the vector $v(\Delta)$ formed by the columns of the upper right triangle of a distance matrix $\Delta$. A certain distance matrix $\Delta$ of a point set $\Ps$ (induced by a labeling of the points) is said to be canonical or maximal if its vector $v(\Delta)$ is the largest one in the set of all vectors of distance matrices of $\Ps$ with respect to the lexicographic order.

\begin{figure}[h]
  \begin{minipage}{6.4cm}
    \begin{center}
      \setlength{\unitlength}{0.0027mm}
      \begin{picture}(23000,19000)
        \qbezier(22270,1000)(11135,1000)(0,1000)
        \qbezier(11731.9,19726.5)(5865.97,10363.27)(0,1000)
        \qbezier(11731.9,19726.5)(17001,10363.27)(22270,1000)
        \qbezier(14433.1,9274.96)(7216.56,5137.48)(0,1000)
        \qbezier(14433.1,9274.96)(18351.6,5137.48)(22270,1000)
        \qbezier(14433.1,9274.96)(13082.5,14500.7)(11731.9,19726.5)
        \qbezier(7910.18,5791.09)(3955.09,3395.55)(0,1000)
        \qbezier(7910.18,5791.09)(15090.1,3395.55)(22270,1000)
        \qbezier(7910.18,5791.09)(9821.06,12758.8)(11731.9,19726.5)
        \qbezier(7910.18,5791.09)(11171.6,7533.02)(14433.1,9274.96)
        \qbezier(3298.12,9274.96)(1649.06,5137.48)(0,1000)
        \qbezier(3298.12,9274.96)(12784.1,5137.48)(22270,1000)
        \qbezier(3298.12,9274.96)(7515.03,14500.7)(11731.9,19726.5)
        \qbezier(3298.12,9274.96)(8865.62,9274.96)(14433.1,9274.96)
        \qbezier(3298.12,9274.96)(5604.15,7533.02)(7910.18,5791.09)
        \qbezier(8567.15,-88.32)(4283.58,455.842)(0,1000)
        \qbezier(8567.15,-88.32)(15418.6,455.842)(22270,1000)
        \qbezier(8567.15,-88.32)(10149.5,9819.11)(11731.9,19726.5)
        \qbezier(8567.15,-88.32)(11500.1,4593.32)(14433.1,9274.96)
        \qbezier(8567.15,-88.32)(8238.66,2851.39)(7910.18,5791.09)
        \qbezier(8567.15,-88.32)(5932.63,4593.32)(3298.12,9274.96)
        \put(0,1000){\circle*{600}}
        \put(22270,1000){\circle*{600}}
        \put(11731.9,19726.5){\circle*{600}}
        \put(14433.1,9274.96){\circle*{600}}
        \put(7910.18,5791.09){\circle*{600}}
        \put(3298.12,9274.96){\circle*{600}}
        \put(8567.15,-88.32){\circle*{600}}
      \end{picture}  
    \end{center}
  \end{minipage}
  \begin{minipage}[c]{4cm}
    \scriptsize
    \begin{eqnarray*}
      &&\Big(\frac{0}{1},\frac{0}{1}\sqrt{2002}\Big)\\
      &&\Big(\frac{22270}{1},\frac{0}{1}\sqrt{2002}\Big)\\
      &&\Big(\frac{26127018}{2227},\frac{932064}{2227}\sqrt{2002}\Big)\\
      &&\Big(\frac{245363}{17},\frac{3144}{17}\sqrt{2002}\Big)\\
      &&\Big(\frac{17615968}{2227},\frac{238464}{2227}\sqrt{2002}\Big)\\
      &&\Big(\frac{56068}{17},\frac{3144}{17}\sqrt{2002}\Big)\\
      &&\Big(\frac{19079044}{2227},-\frac{54168}{2227}\sqrt{2002}\Big)
    \end{eqnarray*}
  \end{minipage}
  \caption{First example of an integral heptagon in general position.}
  \label{fig_ex_1}
\end{figure}
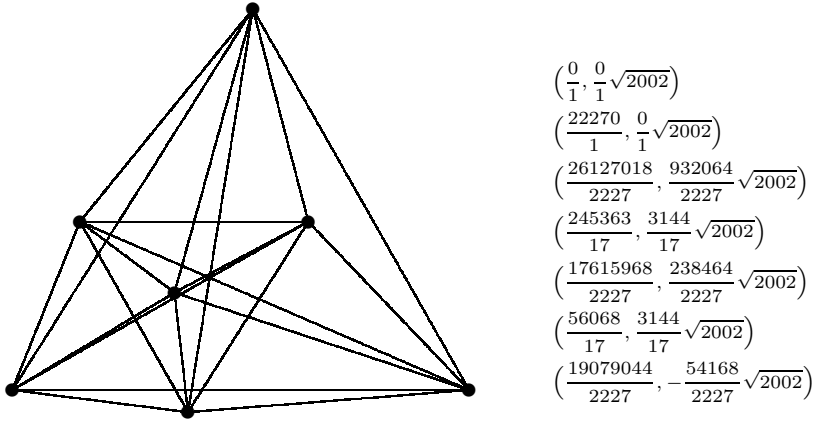

\noindent
In Figure \ref{fig_ex_1} we give an embedding of distance matrix (\ref{example_1}) in the plane and an exact coordinate representation. Discovering this point set clearly motivates to search for further examples to get ideas how to construct an infinite family of examples. Unfortunately this point set is the only example with at most $30000$ in diameter. For diameters greater than $30000$ our approach of exhaustive search requires too much computational power so that we decided to skip to a restricted search. To describe the details of our restriction of the search space we need:

\begin{definition}
  \label{def_characteristic}
  The \textbf{characteristic} of an integral triangle with side lengths $a,b,c\in\mathbb{Z}$ is the square free part of 
  $(a+b+c)(a+b-c)(a-b+c)(-a+b+c)$.
\end{definition}

\begin{theorem}
  Each non degenerated triangle in a plane integral point set has equal characteristic.
\end{theorem}

\noindent
In point set (\ref{example_1}) the characteristic is given by $2002=2\cdot 7\cdot 11\cdot 13$ which explains the shape of the $y$-coordinates, see Figure \ref{fig_ex_1} and \cite{paper_characteristic}. We notice that the characteristic of point set (\ref{example_1}) is composed of relatively small prime factors. By a look at our list of integral hexagons in general position \cite{hp} we see that this seems to be a phenomenon that holds for a great part of the known examples. This phenomenon seems to hold for similar problems also. By determing the minimum diameter $d(2,n)$ of plane integral point sets without further restrictions up to $n=122$ points \cite{paper_alfred} we could check that the known minimal examples also have a characteristic composed of small prime factors. If additionally no three points are allowed to be collinear we denote the corresponding minimum diameter by $\overline{d}(n,2)$. By determing all those minimal integral point sets with up to $n=36$ points \cite{1088.52011,paper_alfred} we could check that the same phenomenon also occurs in this case. So it seemed worth a try to exhaustively construct all plane integral point sets in general position with given diameter of at most $70000$ and the characteristic being a divisor of $6469693230=2\cdot 3\cdot 5\cdot 7\cdot 11\cdot 13\cdot 17\cdot 19\cdot 23\cdot 29$. The outcome was yet another example:

\begin{equation}
  \label{example_2}
  \left(
  \begin{array}{rrrrrrr}
        0 & 66810 & 66555 & 66294 & 49928 & 41238 & 40290 \\
    66810 &     0 & 32385 & 64464 & 32258 & 25908 & 52020 \\
    66555 & 32385 &     0 & 34191 & 16637 & 33147 & 33405 \\
    66294 & 64464 & 34191 &     0 & 34322 & 53244 & 26724 \\
    49928 & 32258 & 16637 & 34322 &     0 & 20066 & 20698 \\
    41238 & 25908 & 33147 & 53244 & 20066 &     0 & 32232 \\
    40290 & 52020 & 33405 & 26724 & 20698 & 32232 &     0
  \end{array}
  \right)
\end{equation}

\noindent
Unfortunately the discovery of further examples is currently beyond our means since the algorithm we use is of running time $\Omega(d^3)$ for the search for plane integral point sets in general position with diameter at most $d$. Though the restriction on the characteristic did accelerate computations significantly the theoretic lower bound for the complexity remains. (There are $O(d^3)$ integral triangles with diameter at most $d$.)

\section{Open problems}

Clearly, one can ask for further examples or an infinite family of integral heptagons in general position. Since our two given examples are in non convex position it would be interesting to see a convex example. As a further restriction Bell and Noll \cite{cluster} also required the coordinates of the point sets to be integral. Such point sets are commonly called $n_m$-clusters, where $n$ is the number of points and $m$ the dimension. In general the set of $n_2$-cluster equals the set of plane integral point sets in general position with characteristic $1$. So far no $7_2$-cluster is known and even its existence is unclear. The smallest $6_2$-cluster has diameter $1886$. At first sight it seems that we have answered Erd\H{o}s question completely, but from a realistic point of view we have only pushed the frontier a step further. Originally P. Erd\H{o}s asked for five points in the plain, no three on a line, no for on a circle with pairwise integral distances. When such a set was found he asked for $6$-set then for a seven set. So now we ask as a substitute:
\begin{center}
  {\glqq}Are there eight points in the plane, no three on a line, no four on a circle with pairwise integral distances?{\grqq}
\end{center}
       
\bibliographystyle{abbrv}
\bibliography{rare}
\bibdata{rare}

\end{document}